\documentclass[a4paper,oneside,reqno]{amsart}

\calclayout

\usepackage{etoolbox}

\newtoggle{default}\toggletrue{default}
\newtoggle{birk}\togglefalse{birk}
\newtoggle{birk_t2}\togglefalse{birk_t2}
\newtoggle{siam}\togglefalse{siam}
\newtoggle{springer}\togglefalse{springer}

\newtoggle{endwithsymbol}\togglefalse{endwithsymbol}

\toggletrue{default}
\toggletrue{endwithsymbol}

\usepackage[utf8]{inputenc}
\usepackage[T1]{fontenc}
\usepackage{lmodern}
\usepackage[english]{babel}
\usepackage{microtype} % better linebreaks

\usepackage{amsmath,amssymb,amsthm}
\usepackage{thmtools} % extends amsthm. for example /autoref (hyperref package required) is very handy

\usepackage[hypertexnames=false]{hyperref} % for internal and external links
\hypersetup{hidelinks}
\usepackage[nameinlink,sort]{cleveref} % for \cref and \Cref (has to be after amsmath)
\usepackage{orcidlink}

\usepackage{mathtools} % many additional math commands
\usepackage{upgreek}
\usepackage{accents} % for \accentset{symbol}
\usepackage{stmaryrd} % for double brackets, e.g., \llbracket \rrbracket

\usepackage{tikz} % for drawing nice pictures
\usetikzlibrary{cd}
\usetikzlibrary{calc}

\usepackage{dsfont} % for indicator function symbol (\mathds{1})
\usepackage{pifont} % for special square (\ding{113} creates the square)

\usepackage{ifthen}
\usepackage{enumitem}
\setlist[itemize]{labelindent=\parindent,leftmargin=*,itemsep=3pt,topsep=5pt}
\setlist[enumerate]{label=\textup{(\roman{*})},itemsep=3pt,labelindent=0pt,leftmargin=*}

\usepackage{nicematrix}

\usepackage{graphicx} % to include external graphics
\usepackage{subcaption}
\usepackage{booktabs} % better tables
\usepackage{esint} % better integral signs

\usepackage{bbm}

\usepackage{color}

\usepackage{pgfpages} % Sebastian
\usepackage{pgfplots} % Sebastian
\pgfplotsset{compat=1.13} % Sebastian

\definecolor{NathanaelColor}{rgb}{0.2,0.0,0.6}

\renewcommand{\leq}{\leqslant}
\renewcommand{\geq}{\geqslant}

\newcommand{\R}{\mathbb{R}} % real numbers
\newcommand{\N}{\mathbb{N}} % natural numbers starting with 1
\newcommand{\Nnull}{\mathbb{N}_{0}} % natural numbers starting with 0
\newcommand{\C}{\mathbb{C}} % complex numbers

\renewcommand{\i}{\mathrm{i}}
\newcommand{\e}{\mathrm{e}}

\newcommand{\Cc}[1][\infty]{\mathrm{C}_{\mathrm{c}}\ifthenelse{\equal{#1}{}}{}{^{#1}}}
\newcommand{\Lp}[2][]{\mathrm{L}_{#2\ifthenelse{\equal{#1}{}}{}{,#1}}} % Lp spaces
\DeclareMathOperator{\rnk}{rnk}

\DeclareMathOperator{\dom}{dom}
\DeclareMathOperator{\diag}{diag}

\DeclarePairedDelimiterX{\dset}[2]{\{}{\}}{#1\,\delimsize\vert\,\mathopen{} #2}
\DeclarePairedDelimiterX{\scprod}[2]{\langle}{\rangle}{#1,#2}

\renewcommand{\Re}{\operatorname{Re}}

\newcommand{\1}{\mathbf{1}}
\renewcommand{\d}{\mathrm{d}}

\theoremstyle{plain}% default
\newtheorem{theorem}{Theorem}[section]

\newtheorem{lemma}[theorem]{Lemma}
\newtheorem{proposition}[theorem]{Proposition}
\newtheorem{corollary}[theorem]{Corollary}

\iftoggle{endwithsymbol}{%
    \declaretheorem[style=definition,sibling=theorem,qed=\ding{169}]{definition}
    \declaretheorem[style=definition,sibling=theorem,qed=\ding{169}]{example}
    \declaretheorem[style=definition,sibling=theorem,qed=\ding{169}]{problem}
    \declaretheorem[style=definition,sibling=theorem,qed=\ding{169}]{assumption}
    \declaretheorem[style=definition,numbered=no,qed=\ding{169}]{claim}

    \declaretheorem[style=remark,sibling=theorem,qed=\ding{169}]{remark}
  }%
  {%
    \theoremstyle{definition}

    \theoremstyle{remark}
    \newtheorem{remark}[theorem]{Remark}
  }%

\allowdisplaybreaks

\begin{document}

% ----------------------------------------------------------------------
% Titlepage
% ----------------------------------------------------------------------
\title[exponential stability for pHS]{Exponential stability of port-Hamiltonian systems with unstructured, measurable Hamiltonian densities} %TODO find good title

\author[H.~Kuettler]{Heinrich Kuettler}
\address[H.K.]{ Palo Alto\\ CA\\ USA}
\email{kuettler@math.lmu.de}

\author[M.~Waurick]{Marcus Waurick\,\orcidlink{0000-0003-4498-3574}}
\address[M.W.]{TU Bergakademie Freiberg \\
  Institute of Applied Analysis \\
  Akademiestrasse 6 \\
  D-09596 Freiberg \\
  Germany \\ phone: +49 3731 39-2688\\ \emph{Corresponding Author}}
\email{marcus.waurick@math.tu-freiberg.de}

%\email{marcus.waurick@math.tu-freiberg.de}

\date{\today}
\dedicatory{}


\ifboolexpr{togl{birk_t2} or togl{birk}}{%
}{%
}%

\thanks{}

% ----------------------------------------------------------------------

\begin{abstract} For first-order port-Hamiltonian systems in the dissipative setting considered here, we characterise exponential stability without any further regularity or block-structure assumption on the Hamiltonian density, thus closing a gap persisting in the literature. The proof requires revisiting the computation of the adjoint generator of the port-Hamiltonian semigroup. Finally, we provide both examples and non-examples showcasing the inconclusiveness of known conditions on the Hamiltonian densities necessary for the application of characterisation results in the literature.
 \end{abstract}

\ifboolexpr{togl{default} or togl{birk_t2} or togl{birk}}{\maketitle}{}%

\section{Introduction}\label{sec:intro}

In recent years, stability properties for example systems have shown growing interest. In particular, due to the possibility of invoking arguments from linear algebra, a particular focus was given to a class of $1$-D hyperbolic partial differential equations often referred to as port-Hamiltonian systems -- or sometimes called distributed port-Hamiltonian systems\footnote{One of the mathematical sources for these systems is \cite{JZ12}, where these systems were called `port-Hamiltonian systems'. Hence, we shall employ the same custom here; see also the discussion in \cite{SW26}.}. We refer to \cite{v06,vJ14,JZ12,S21} for some classical and more recent sources on port-Hamiltonian systems. In the literature, one can now find a characterisation of strong stability for these systems (\cite{WZ24}) and a characterisation of exponential stability (\cite{TW26}); quantified strong stability results are available in \cite{ASVW26}. Whereas the former criterion invokes only the invertibility of a family of matrices, the latter two require the matrix norm of the inverses to satisfy certain bounds across the entire family. 

The characterisation from \cite{WZ24} applies to all port-Hamiltonian systems in the setting considered here. In contrast, the ones in \cite{TW26,ASVW26} ask for a regularity condition satisfied by the port-Hamiltonian system in question. More precisely, the asked-for property is the boundedness of an associated family of fundamental solutions, which holds, for instance, if the Hamiltonian density is of bounded variation or leaves the relevant spectral subspaces invariant. In any case, the regularity condition has been subject to an open problem in \cite[Problem 1.6]{TW26}, which asks as to whether this condition is needed to obtain exponential stability results in the first place. We will answer this question in the negative. In fact, more profoundly, the present article resolves the issue of the characterisation of exponential stability for port-Hamiltonian systems irrespective of the regularity property. Hence, we provide the full picture of the characterisation of exponential stability, closing the gap in the literature. It remains to analyse the quantified strongly stable decay of systems violating the regularity condition.  

For the characterisation result of exponential stability here, the strategy is to provide a precise description of the semigroup generator of port-Hamiltonian systems, which will particularly facilitate finding the domain of the adjoint of the generator -- in contrast to earlier findings (\cite[Proposition 3.4.3]{A16} and \cite[Theorem 2.24]{V07}), the present computation of the adjoint is particularly useful here as it avoids the introduction of effort- and flow-variables and clearly identifies the dimension of `dissipating channels'. The latter number is the same for the generator and its adjoint. In particular, they share a common kernel at the boundary, which may be interpreted as a periodicity-like boundary condition on the non-strictly dissipative part. 

With this information, we are in a position to employ some theory of dissipative operators, more precisely, for functional models as provided in \cite[Section 4]{BMNW20}. A distance estimate of paths of unitary operators (see also \cite{ALV14}) will then be decisive for the resolvent estimate to eventually invoke the Gearhart--Pr\"uss theorem (see \cite{P84}) establishing the characterisation for exponential stability.

In \Cref{sec:pHsmr}, we provide an overview of the functional analytic set-up of port-Hamiltonian systems and recall the standard generation theorem. In particular, we revisit the integration-by-parts formula, which we will employ when discussing the adjoint generator. Before the computation of the adjoint, we need some results in linear algebra to be presented in \Cref{sec:sigmadiss}. The first non-trivial result of the present note is the formula for the adjoint of the port-Hamiltonian generator in \Cref{sec:adjgen}. \Cref{sec:charexp} describes the existing results on exponential stability and the new contribution and result of this manuscript. The subsequent section, \Cref{sec:res}, contains the proof of the main result and some applications for special cases. In addition to the characterisation of exponential stability without the regularity condition, we also provide a new characterisation of exponential stability similar to the one in \cite{TW26} -- with a weaker hypothesis. We conclude by providing examples in \Cref{sec:ex}.

%The next section will provide the precise mathematical notions and main results of the paper. Moreover, we will recall some results from the literature to provide the context and detail the contribution of the present text. \Cref{sec:ipb} is focused on the integration by parts formula and a new description of the domain of the port-Hamiltonian semigroup generator and, as a consequence, of its adjoint. The proof of the main results and the necessary estimates for these can be found in \Cref{sec:res}. Examples are treated in \Cref{sec:ex}.

\section{port-Hamiltonian systems}\label{sec:pHsmr}

We recall the setting of port-Hamiltonian systems provided in \cite{JZ12} next. Let $d\in \N$, $P_0=-P_0^*,P_1=P_1^* \in \C^{d\times d}$, $P_1$ invertible, $\mathcal{H}\in L_\infty(a,b)^{d\times d}$ satisfying $\mathcal{H}(x)=\mathcal{H}(x)^*\geq c$ for almost every $x\in (a,b)$ and some $c>0$ in the sense of positive definiteness. Introduce the Hilbert space
\[
  X\coloneqq  L_{2,\mathcal{H}}(a,b)^d \text{ with }\langle u,v\rangle_X\coloneqq \langle \mathcal{H}u,v\rangle_{L_2(a,b)^d}.
\] 
Throughout this article, we will use
\begin{align*}
   A_{\textnormal{max}} \colon \dom(   A_{\textnormal{max}})\subseteq X &\to X \\
   u& \mapsto P_1 (\mathcal{H} u)'+P_0 \mathcal{H} u,
\end{align*}
where
\[
   \dom(A_{\textnormal{max}} ) = \{ u\in X; (\mathcal{H}u)' \in X\},
\]
with the derivative understood in the distributional sense. It is not difficult to see that $A_{\textnormal{max}}$ is closed. Thus, endowed with the graph norm, $   \dom(A_{\textnormal{max}} )$ becomes a Hilbert space in its own right. By the Rellich--Kondrachov theorem,
\[
   \dom(A_{\max})\hookrightarrow X\text{ compactly}.
\]Due to the one-dimensional Sobolev embedding theorem, the mapping 
\[
\dom(A_{\textnormal{max}}) \to C[a,b]^d, u\mapsto \mathcal{H}u,
\] of choosing a continuous representative of $\mathcal{H}u$, is well-defined and continuous. As a consequence, the mapping
\[
   \gamma_{\mathcal{H}} \colon  \dom(A_{\textnormal{max}} ) \to \C^{2d}, u \mapsto \begin{pmatrix}(\mathcal{H}u)(b) \\ (\mathcal{H}u)(a)\end{pmatrix}
\]
is well-defined. If $\mathcal{H}$ is clear from the context we shall write $\gamma$ instead of $\gamma_{\mathcal{H}}$.
The starting point in the analysis of port-Hamiltonian systems is a generation result due to \cite{ZGMV10}. This result provides a neat description of all dissipative restrictions of $A_{\textnormal{max}}$. 
The generation theorem is based on an integration-by-parts formula, which is itself also a consequence of the Sobolev embedding theorem.
We will need
\[
\Sigma \coloneqq \begin{pmatrix} P_1 & 0 \\ 0 & -P_1\end{pmatrix}.
\]
\begin{proposition}\label{prop:ibp} Let $u,v\in \dom(A_{\max})$. Then
\[
    \langle A_{\max} u, v\rangle_X +     \langle u, A_{\max} v\rangle_X = \langle 
    \Sigma \gamma u, \gamma v\rangle_{\mathbb{C}^{2d}}.
\]
\end{proposition}
\begin{proof}
By skew-selfadjointness of $P_0$, we get, using the  standard integration by parts formula together with $P_1=P_1^*\in \C^{d\times d}$,
\begin{align*}
    \langle A_{\max} u, v\rangle_X +     \langle u, A_{\max} v\rangle_X & =
    \langle P_1 (\mathcal{H}u)', \mathcal{H}v\rangle +     \langle \mathcal{H}u, P_1(\mathcal{H}v)'\rangle \\
    & = \langle P_1 (\mathcal{H}u)(b), (\mathcal{H}v)(b)\rangle_{\mathbb{C}^d} - \langle P_1 (\mathcal{H}u)(a), (\mathcal{H}v)(a)\rangle_{\mathbb{C}^d} 
\\
& = \langle 
    \Sigma \gamma u, \gamma v\rangle_{\mathbb{C}^{2d}}.     \qedhere
\end{align*}
\end{proof}
The generation theorem now reads as follows.
\begin{theorem}[\cite{ZGMV10,JZ12}]\label{thm:gen}  Let $W \in \C^{d\times 2d}$ with $\rnk W=d$. The operator $A\subseteq A_{\textnormal{max}}$ given by 
\[
   \dom(A) = \{ u\in \dom(A_{\textnormal{max}}); \gamma_{\mathcal{H}} u \in \ker(W)\}
\]
is m-dissipative on $X$ (and, thus, generates a $C_0$-semigroup of contractions), if, and only if,
\begin{equation}\label{eq:sig}
\forall q\in \ker(W)\colon \Re \langle \Sigma q,q\rangle \leq 0.
\end{equation}
\end{theorem}
\begin{proof}
By \cite[Theorem 7.2.4]{JZ12}, $A$ is m-dissipative if and only if for all $u\in \dom(A)$ we have
\begin{equation}\label{eq:sigreal}
 \Re \langle \Sigma \gamma u,\gamma u\rangle =  2\Re \langle Au,u\rangle_X\leq 0,
\end{equation}
where we used \Cref{prop:ibp} in the first equality. It is not difficult to see that $\dom(A_{\max})\ni u \mapsto \gamma u \in \C^{2d}$ is onto. Hence, $\gamma|_{\dom(A)}$ maps onto $\ker(W)$. 

Now, m-dissipativity of $A$ together with \Cref{eq:sigreal} and the mentioned surjectivity property leads to \Cref{eq:sig}. On the other hand, \Cref{eq:sig} in conjunction with \Cref{eq:sigreal} leads to $\Re \langle Au,u\rangle \leq 0$ for $u\in \dom(A)$, which is sufficient for $A$ to be m-dissipative.
\end{proof}

\section{$\Sigma$-dissipative subspaces}\label{sec:sigmadiss}

In \Cref{thm:gen}, m-dissipative restrictions of the maximal port-Hamiltonian operator are characterised by the boundary values belonging to a $\Sigma$-dissipative, $d$-dimensional subspace, where $\mathcal{W}\subseteq \C^{2d}$ is called \textbf{$\Sigma$-dissipative} if
\[
   \forall q\in \mathcal{W}\colon \Re \langle \Sigma q,q\rangle \leq 0,
\]
with $\Sigma=\diag(P_1,-P_1)$ with $P_1=P_1^*$ self-adjoint and invertible, as before. This section is devoted to a description of $\Sigma$-dissipative subspaces. The main observation reads as follows.
\begin{theorem}\label{thm:sigdiss} Let $\mathcal{W}\subseteq \C^{2d}$ be a $d$-dimensional, $\Sigma$-dissipative subspace. Then there exists an invertible $Z\in \C^{2d\times 2d}$, $r\in \{0,\ldots,d\}$, $\Lambda_1 =\diag(\lambda_1,\ldots,\lambda_r)\in [0,1)^{r\times r}$ such that
\[
   q \in \mathcal{W} \iff  \begin{pmatrix} 1_d & - \begin{pmatrix} 1_{d-r} & 0 \\ 0 & \Lambda_1\end{pmatrix}\end{pmatrix} Zq=0.
\]
Moreover, writing  $Zq=(Z_+q,Z_-q)\in \C^{d}\times \C^d$ and $Z_{\pm}q = (\check{Z}_{\pm}q,\hat{Z}_{\pm}q) \in \C^{d-r}\times \C^r$ we have
\begin{align*}
 \Re \langle \Sigma q ,q\rangle & = \Re \langle \sigma(\check{Z}_+q,\check{Z}_-q),\delta  (\check{Z}_+q,\check{Z}_-q)\rangle \\ & \quad +  \|(1-\Lambda_1^2)^{-1/2}(\hat{Z}_+q- \Lambda_1 \hat{Z}_-q)\|^2  - \|(1-\Lambda_1^2)^{-1/2}(\hat{Z}_-q - \Lambda_1 \hat{Z}_+q)\|^2,
\end{align*}
where $\sigma(\xi,\eta)=\xi+\eta$ and $\delta(\xi,\eta)=\xi-\eta$. Moreover,
\[
   \langle \Sigma q, v\rangle = \langle Z_{+}q,Z_{+}v\rangle -\langle Z_{-}q,Z_{-}v\rangle\quad(q,v\in \C^{2d}).
\]
\end{theorem}

\begin{proof} Let $\tilde{\Sigma}_+$ and $\tilde{\Sigma}_-$ be the positive and negative spectral parts of $\Sigma$, respectively. Then $\Sigma=\tilde{\Sigma}_+ -\tilde{\Sigma}_-$. We identify each of its $d$-dimensional positive and negative spectral subspaces unitarily with $\C^d$. Moreover, with $q_{\pm}\coloneqq \Sigma_{\pm}q\coloneqq \sqrt{\tilde{\Sigma}}_{\pm} q$, we obtain
\[
    \langle \Sigma q ,q\rangle = \|q_+\|^2 - \|q_-\|^2. 
\]
Hence,
\begin{equation}\label{eq:qpqm}
\forall q\in \mathcal{W}\colon   \|q_+\| \leq \|q_-\|.
\end{equation}
Since $P_1$ is invertible, the mapping $q\mapsto (q_+,q_-)$ is invertible. Moreover, the latter inequality shows that $Q_-\colon \mathcal{W}\ni q\mapsto q_- \in \C^d$ is injective. Indeed, for $q\in \mathcal{W}$, $q_-=0$ implies $q_+=q_-=0$, which in turn yields $q=0$. As $\mathcal{W}$ is $d$-dimensional, $Q_-$ is onto and, thus, bijective. Hence, 
\[
   M \colon \C^d \to \C^d, w \mapsto \Sigma_+ Q_-^{-1} w
\]is well-defined. We have $\|M\|\leq 1$. Indeed, let $w\in \C^d$. Then we find a unique $q \in \mathcal{W}$ such that $q = Q_-^{-1} w$ and, thus, $q_-=w$. As 
\[
Mq_- = Mw=\Sigma_+ Q_{-}^{-1} w=\Sigma_+ q = q_+,
\] contractivity of $M$ follows from \Cref{eq:qpqm}.

Next, for $q\in \C^{2d}$,
\begin{equation}\label{eq:spsm}
   q \in \mathcal{W} \iff \Sigma_+q - M\Sigma_-q=0.
\end{equation}
Indeed, if $q\in \mathcal{W}$, then $Q_-^{-1}\Sigma_-q=q$ and, thus, $M\Sigma_- q = \Sigma_+ Q_-^{-1} \Sigma_- q =q_+$. If $\Sigma_+q - M\Sigma_-q=0$, then $\Sigma_+ q = \Sigma_+  Q_-^{-1} \Sigma_- q $. In particular, $\tilde{q}\coloneqq Q_-^{-1} \Sigma_- q \in \mathcal{W}$ satisfies $\Sigma_-\tilde{q} =q_-$ and $\Sigma_+ \tilde{q} =\Sigma_+ q$. Since $q\mapsto (q_+,q_-)$ is bijective, it follows that $\tilde{q}=q$. Hence, \Cref{eq:spsm} is established.

Since $\|M\|\leq 1$, we obtain $1-M^*M\geq 0$ in the sense of positive definiteness. Hence, we now find a unitary matrix $T\in \C^{d\times d}$, $r\in \{0,\ldots,d\}$ and a real diagonal matrix $\Lambda=\diag(1_{d-r},\Lambda_1)$  with $\Lambda_1=\diag(\lambda_1,\ldots, \lambda_{r}) \in \C^{r\times r}$ with $0\leq \lambda_j<1$ for all $j\in \{1,\ldots, r\}$ such that $M^*M=T^*\Lambda^2 T$. Then
\[
    \langle (1-M^*M)q_-,q_-\rangle =    \langle (1-\Lambda^2)Tq_-,Tq_-\rangle\quad(q\in \mathcal{W}).
\]
Next, let $q\in \C^{2d}$. Taking a unitary $U$ such that $M=U|M|=U\sqrt{M^*M}$, we put  $Z_-q\coloneqq v_- \coloneqq Tq_{-}$ and  $Z_+q\coloneqq v_+ \coloneqq TU^*q_{+}$ and $\check{v}_\pm \coloneqq  \big(v_{\pm,j}\big)_{j\in \{1,\ldots,d-r\}}$ and $\hat{v}_\pm \coloneqq \big(v_{\pm,j}\big)_{j\in \{d-r+1,\ldots,d\}}$. Then $Z$ is invertible. We establish the reformulation of $\Re \langle \Sigma q,q\rangle$ first. For this, since $\Lambda_1$ is diagonal, the formula
\begin{align*}
\|x-\lambda y\|_{\C}^2-\|y-\lambda x\|_{\C}^2 = \langle (\sqrt{1-\lambda^2}) x, (\sqrt{1-\lambda^2}) x\rangle_\C - \langle (\sqrt{1-\lambda^2}) y, (\sqrt{1-\lambda^2})y\rangle_\C
\end{align*}
for  $x,y \in \C$ and $\lambda \in [0,1]$, implies 
\[
 \|\hat{v}_+\|^2 - \|\hat{v}_-\|^2 =\|(1-\Lambda_1^2)^{-1/2}(\hat{v}_+ - \Lambda_1 \hat{v}_-)\|^2  - \|(1-\Lambda_1^2)^{-1/2}(\hat{v}_- - \Lambda_1 \hat{v}_+)\|^2.
\]
Moreover, recalling  $\delta  (u,v)\coloneqq u-v$ as well as $\sigma(u,v) \coloneqq u+v$, we infer
\[
    \Re \langle \sigma(u,v),\delta  (u,v)\rangle = \langle  u,u\rangle - \langle  v,v\rangle.
    \]
  Thus,
    \begin{align*}
    \Re \langle \Sigma q,q\rangle 
        & = \|q_+\|^2 - \|q_-\|^2\\
        &= \|TU^*q_+\|^2-\|Tq_-\|^2 = \|v_+\|^2 - \|v_-\|^2  \\
        &= \|\check{v}_+\|^2 - \|\check{v}_-\|^2 + \|\hat{v}_+\|^2 - \|\hat{v}_-\|^2  
         \\ &=\Re \langle \sigma(\check{v}_+,\check{v}_-),\delta  (\check{v}_+,\check{v}_-)\rangle \\ & \quad +  \|(1-\Lambda_1^2)^{-1/2}(\hat{v}_+ - \Lambda_1 \hat{v}_-)\|^2  - \|(1-\Lambda_1^2)^{-1/2}(\hat{v}_- - \Lambda_1 \hat{v}_+)\|^2,
    \end{align*}
which establishes the claimed formula. The formula for the scalar product follows from polarisation and 
\[
 \Re \langle \Sigma q,q\rangle = \|TU^*q_+\|^2-\|Tq_-\|^2 = \|Z_+q\|^2-\|Z_-q\|^2.
\] Finally, we are concerned with the equivalence for which we invoke \Cref{eq:spsm}. Thus, we need to show
    \[
       q \in \mathcal{W} \iff \hat{v}_+ - \Lambda_1 \hat{v}_-=0 \text{ and }\check{v}_+=\check{v}_-.
    \]
    The conditions
    \[
    \hat{v}_+ - \Lambda_1 \hat{v}_-=0 \text{ and }\check{v}_+=\check{v}_-
    \]
    are the same as saying
    \[
      v_+  = \Lambda v_-,\text{ i.e., } TU^*q_+ = \Lambda Tq_-\text{ or }U^*q_+ = T^*\Lambda Tq_- = \sqrt{M^*M} q_- = |M|q_-.
    \]
    The latter is the same as saying $q_+=Mq_-$.\end{proof}
    \begin{remark}\label{rem:z-onto}
    Another look at the proof of \Cref{thm:sigdiss} reveals that $\mathcal{W}\ni q\mapsto Z_-q$ is onto.
    \end{remark}
    \begin{corollary}\label{cor:sigdissadj} Let $\mathcal{W}\subseteq \C^{2d}$ be a $d$-dimensional, $\Sigma$-dissipative subspace and let $r,\Lambda_1, Z$ be as in \Cref{thm:sigdiss}. Let $v\in \C^{2d}$. Then 
    \[
       \langle \Sigma q, v\rangle =0 \quad (q\in \mathcal{W}) \iff \check{Z}_+v=\check{Z}_-v\text{ and }\Lambda_1\hat{Z}_+v=\hat{Z}_-v.
       \]
    \end{corollary}
  \begin{proof}
    By \Cref{thm:sigdiss}, we have 
    \[
    q\in \mathcal{W} \iff \check{Z}_+q=\check{Z}_-q \text{ and }\hat{Z}_+q=\Lambda_1\hat{Z}_-q.
    \]
    Thus, let $v\in \C^{2d}$ be such that 
    \[
        \langle \Sigma q, v\rangle =0 \quad (q\in \mathcal{W}).
    \]
    Then, using \Cref{thm:sigdiss}, we get
    \begin{align*}
   &   \langle \Sigma q, v\rangle  = \langle Z_+q,Z_+v\rangle- \langle Z_-q,Z_-v\rangle\\
     & = \langle \check{Z}_+q,\check{Z}_+v\rangle- \langle \check{Z}_-q,\check{Z}_-v\rangle + \langle \hat{Z}_+q,\hat{Z}_+v\rangle- \langle \hat{Z}_-q,\hat{Z}_-v\rangle \\
      & = \langle \check{Z}_-q,\check{Z}_+v\rangle- \langle \check{Z}_-q,\check{Z}_-v\rangle + \langle \Lambda_1\hat{Z}_-q,\hat{Z}_+v\rangle- \langle \hat{Z}_-q,\hat{Z}_-v\rangle \\
            & = \langle \check{Z}_-q,\big(\check{Z}_+v- \check{Z}_-v\big)\rangle + \langle \hat{Z}_-q,\big(\Lambda_1\hat{Z}_+v-\hat{Z}_-v\big)\rangle 
    \end{align*}
    Since $\mathcal{W}\ni q\mapsto Z_-q$ is onto by \Cref{rem:z-onto}, the claim follows (the computation in reverse direction establishes the equivalence).
    \end{proof}

\section{The adjoint generator}\label{sec:adjgen}

As a first application of the linear algebra result from the previous section, we may now compute the adjoint of the port-Hamiltonian operator. We adopt the setting in \Cref{sec:pHsmr}. Throughout this section, we let $W\in \C^{d\times 2d}$ with $\rnk W=d$. Moreover, we let $A\subseteq A_{\max}$ such that
\[
   \dom(A) = \{u\in \dom(A_{\max}); \gamma u \in \ker W\}.
\]
Then a reformulation of \Cref{thm:gen} reads as follows.
\begin{theorem}\label{thm:gen2} $A$ generates a $C_0$-semigroup of contractions, if, and only if, $\ker W$ is $\Sigma$-dissipative.
\end{theorem}
Thus, henceforth, we shall assume that $\ker W$ is $\Sigma$-dissipative. Since $W$ has full rank, $\ker W$ is $d$-dimensional. Hence, \Cref{thm:sigdiss} applies and we find $r,\Lambda_1,Z$ such that for $q\in \C^{2d}$
\[
q\in \ker W \iff \check{Z}_+q=\check{Z}_-q \text{ and }\hat{Z}_+q=\Lambda_1\hat{Z}_-q.
\] Hence, we obtain
\[
\dom(A) = \{ u \in \dom(A_{\max}); \check{Z}_+\gamma u=\check{Z}_-\gamma u \text{ and }\hat{Z}_+\gamma u=\Lambda_1\hat{Z}_-\gamma u\}.
\]
The formula for the adjoint of $A$ is similar as we shall see next.
\begin{theorem}\label{thm:adjform} Within the setting of this section, we have $A^*\subseteq -A_{\max}$ and
\[
   \dom(A^*)=\{ v \in \dom(A_{\max}); \check{Z}_+\gamma v=\check{Z}_-\gamma v \text{ and }\Lambda_1\hat{Z}_+\gamma v=\hat{Z}_-\gamma v\}.
\]
\end{theorem}
\begin{proof} It follows directly from skew-selfadjointness of $P_0$, selfadjointness of $P_1$ and the skew-symmetry of $\partial$ restricted to compactly supported test functions, that $A^*\subseteq -A_{\max}$. Next, recall that $\dom(A)\ni u\mapsto \gamma u\in \ker(W)$ is onto. Thus, for $v\in\dom(A_{\max})$, we have $v\in \dom(A^*)$ if and only if for all $u\in \dom(A)$
\[
 0 =   \langle A_{\max} u, v\rangle_X +     \langle u, A_{\max} v\rangle_X \big(= \langle 
    \Sigma \gamma u, \gamma v\rangle_{\mathbb{C}^{2d}}\Big),
\]where the expression in brackets follows from \Cref{prop:ibp}. Thus, 
\[
v\in \dom(A^*)\iff \forall q\in \ker(W) \colon \langle \Sigma q, \gamma v\rangle=0.
\]
Hence, by \Cref{cor:sigdissadj}, the last condition is equivalent to
\[
\check{Z}_+\gamma v=\check{Z}_-\gamma v \text{ and }\Lambda_1\hat{Z}_+\gamma v=\hat{Z}_-\gamma v.\qedhere
\]
\end{proof}
We will write the domains of $A$ and $A^*$ in a more compact form. For this, we introduce
\begin{align}\label{eq:LGH}
\begin{aligned}
  L & \colon \C^{2d} \to \C^{d-r}, q\mapsto \check{Z}_+q-\check{Z}_-q, & L_\gamma& \coloneqq L\circ \gamma \\
  G &\colon \C^{2d} \to \C^{r}, q\mapsto (1-\Lambda_1^2)^{-1/2}(\hat{Z}_+q-\Lambda_1\hat{Z}_-q), & G_\gamma &\coloneqq G\circ \gamma, \\
  H &\colon \C^{2d} \to \C^{r}, q\mapsto (1-\Lambda_1^2)^{-1/2}(\Lambda_1\hat{Z}_+q-\hat{Z}_-q), & H_\gamma& \coloneqq H\circ \gamma.
  \end{aligned}
\end{align} All mappings with the subscript $\gamma$ have domain $\dom(A_{\max})$.
By \Cref{thm:sigdiss} and polarisation, for $q,v\in\ker L$ we have
\[
 \langle\Sigma q,v\rangle=\langle Gq,Gv\rangle-\langle Hq,Hv\rangle.
\]
Then we have
\begin{align*}
  \dom(A) & = \{ u\in \dom(A_{\max}); L_\gamma u = 0, G_\gamma u=0\}=\ker(L_\gamma)\cap \ker(G_\gamma),\\
  \dom(A^*) & = \{ u\in \dom(A_{\max}); L_\gamma u = 0, H_\gamma u=0\}=\ker(L_\gamma)\cap \ker(H_\gamma).
\end{align*}
There is a nice characterisation of the situation when $r=d$:
\begin{corollary} The following conditions are equivalent:
\begin{enumerate}
  \item[(i)] for all $u\in \dom(A)$ with $\gamma u\neq 0$, we have $\Re \langle Au,u\rangle_X<0$;
 \item[(ii)] $r=d$.
\end{enumerate}
\end{corollary}
\begin{proof}From \Cref{thm:sigdiss} and \Cref{prop:ibp}, we deduce for all $u\in \dom(A)$ since $L_\gamma u =0$
\begin{equation}\label{eq:raau}
 2\Re \langle Au,u\rangle = \Re \langle \Sigma \gamma u,\gamma u\rangle = -\| H_\gamma u \|^2.
\end{equation}
For $q\in \ker W$, we have 
\[
 Hq=-(1-\Lambda_1^2)^{1/2}\hat Z_-q.
\]
By \Cref{rem:z-onto}, $Z_-|_{\ker W}$ is a bijection onto
$\C^d$, so $\dim\ker(H|_{\ker W})=d-r$.
Since $\gamma[\dom(A)]=\ker W$, \Cref{eq:raau} shows
that (i) holds if and only if $r=d$.
\end{proof}

\section{Characterisation of Exponential Stability}\label{sec:charexp}
Before we begin with the presentation of the new contributions in the present article concerning exponential stability, we recall some results obtained in the literature. We adopt the setting from \Cref{sec:pHsmr} and let $W\in \C^{d\times 2d}$ with $\rnk W=d$ be such that $\ker W$ is $\Sigma$-dissipative. We will use the formulas in \Cref{eq:LGH}, which particularly includes $r$, $\Lambda_1$, and $Z$.

For the stability theorems we furthermore introduce for $t\in \R$ the fundamental solution $\Phi_t \in C[a,b]^{d\times d}$ satisfying $\Phi_t(a)=1_d$ associated to the ordinary differential equation
\[
   v'(x) =  P_1^{-1}(\i t \mathcal{H}(x)^{-1}  - P_0) v(x)\quad (x\in (a,b)).
\]
With this, we define for $t\in \R$
\[
    T_t \coloneqq W \begin{pmatrix} \Phi_t(b) \\ 1_d\end{pmatrix}.
\]
There is a direct relation of the invertibility of $T_t$ to the resolvent of $A$:
\begin{proposition}[\cite{TW26,WZ24}]\label{prop:TWZ} Let $t\in \R$. Then
\[
    T_t \text{ invertible }\iff \i t \in \rho(A).
\]
\end{proposition}
We denote by $(U_t)_{t\geq 0}$, the semigroup generated by $A$.
\begin{theorem}[\cite{WZ24}]\label{thm:WZ} The following conditions are equivalent:
\begin{enumerate}
\item[(i)] $(U_t)_{t\geq0}$ is semi-uniformly stable;
\item[(ii)] $(U_t)_{t\geq0}$ is strongly stable;
\item[(iii)] for all $t\in \R$, $T_t$ is invertible.
\end{enumerate} 
\end{theorem}
The characterisations in \cite{TW26,ASVW26} require the following additional condition
\begin{equation}\label{eq:(B)}\tag{B}
\sup_{t\in \R} \|\Phi_t\|_{\infty}<\infty.
\end{equation}
The available criterion for exponential stability now reads as follows.
\begin{theorem}[{{\cite[Theorem 3.5]{TW26}}}]\label{thm:WTB} Assume that \eqref{eq:(B)} holds.
Then
\[
 (U_t)_{t\geq 0} \text{ exponential stability }\iff T_t \text{ invertible and }\sup_{t\in \R} \|T_t^{-1}\|<\infty.
\]
\end{theorem}
To present the criterion for exponential stability without \eqref{eq:(B)}, we define
\[
  \hat{W} \colon \C^{2d} \ni q \mapsto (Lq,Gq)\in \C^{d-r}\times \C^r
\]
and note, by the compact description of $\dom(A)$ and the definition of $A$ that
\[
   \ker W = \ker{\hat{W}}.
\]
Thus, we find an invertible $Y$ such that
\[
  \hat{W} = YW.
  \]
  With this, we put\[
   \hat{T}_t \coloneqq Y T_t = YW\begin{pmatrix} \Phi_t(b) \\ 1_d\end{pmatrix}.
\] Introducing 
\[
    \Psi_t \coloneqq \int_a^b \Phi_t(x)^*\mathcal{H}(x)^{-1}\Phi_t(x)\d x\text{ and } J_r\coloneqq \begin{pmatrix} 0_{d-r\times r}\\ 1_r\end{pmatrix},
\]
we may now provide the main theorem of the present article.
\begin{theorem}[Characterisation of exponential stability without \eqref{eq:(B)}]\label{thm:mr-raw} Let $A$ be as above. Then the following conditions are equivalent:
\begin{enumerate}
\item[(i)] $(U_t)_{t\geq 0}$ is exponentially stable;
\item[(ii)]$T_t$ is invertible for all $t\in \R$ and 
\[
  \sup_{t\in \R} \|\Psi_t^{1/2} \hat{T}_t^{-1} J_r \|<\infty,
\]
\end{enumerate}
\end{theorem}
\begin{remark} Let us provide the perspective to \Cref{thm:WTB}. Thus, assume that \eqref{eq:(B)}  holds. Then $\sup_{t\in \R} \|T_t^{-1}\|<\infty$ leads to (ii) in \Cref{thm:mr-raw}, which yields exponential stability. Yet, exponential stability of $(U_t)_{t\geq 0}$ in \Cref{thm:mr-raw} even with \eqref{eq:(B)}  and the help of \cite[Lemma 3.2]{TW26}, directly provides only
\[
   \sup_{t\in \R} \|{T}_t^{-1} Y^{-1}J_r \|<\infty.
\]Yet, by \Cref{thm:WTB}, this turns out to be equivalent to
\[
  \sup_{t\in \R} \|T_t^{-1}\|<\infty.
\]
Thus, under \eqref{eq:(B)}, the two characterisations are equivalent, although the estimate in \Cref{thm:mr-raw} does not yield full control of $T_t^{-1}$ without invoking \Cref{thm:WTB}. For a different relation, we also refer to \Cref{thm:charBpr} below.
\end{remark}
  \section{Proof of the main result}\label{sec:res}
This section is devoted to a proof of \Cref{thm:mr-raw}. For this, we adopt the setting in \Cref{sec:charexp} together with the entire notation of this section. In particular, for $t\in \R$,
\[
    \Phi_t \colon [a,b]\to \C^{d\times d}
\]
is the fundamental solution of 
\[
   v' = P_1^{-1} (\i t \mathcal{H}(\cdot)^{-1}-P_0)v
\]
normalised by $\Phi_t(a)=1_{d}$ and
\[
    T_t \coloneqq W\begin{pmatrix} \Phi_t(b) \\ 1_d\end{pmatrix} \text{ as well as }   \hat{T}_t \coloneqq  \hat{W}\begin{pmatrix} \Phi_t(b) \\ 1_d\end{pmatrix}=Y{W}\begin{pmatrix} \Phi_t(b) \\ 1_d\end{pmatrix}
\]
for some invertible $Y$. At first we take a look at the eigenvalue problem for $A_{\max}$ and its relation to $\Phi_t$.
\begin{lemma}\label{lem:funmat} Let $v_0\in \C^d$, $t\in \R$, $v\coloneqq \mathcal{H}(\cdot)^{-1}\Phi_t(\cdot)v_0$. Then $v\in \dom(A_{\max})$ and 
\[
A_{\max} v = \i t v\text{ and }\gamma v = \begin{pmatrix} \Phi_t(b) \\ 1_d\end{pmatrix} v_0.
\] 
Moreover, if $T_t$ is invertible and  $v_0=\hat{T}_t^{-1}w_0= T_t^{-1} Y^{-1}w_0 \in \C^d$ for some $w_0\in \C^d$, then
\[
 \hat{W} \gamma v =  w_0.
\]
\end{lemma}
\begin{proof}
The first statement follows directly from the properties of the chosen fundamental solution. For the second statement, we use the first one and compute
\[
  \hat{W} \gamma v = YW \begin{pmatrix} \Phi_t(b) \\ 1_d\end{pmatrix} v_0=Y T_t v_0= YT_tT_t^{-1}Y^{-1}w_0 =w_0.\qedhere
\]
\end{proof}
Next, recall $ J_r= \begin{pmatrix} 0_{d-r\times r}\\ 1_r\end{pmatrix}$ and define for $t\in \R$ with $\i t\in\rho(A)$
\[
   R_t \colon \C^r \to X,  z \mapsto \mathcal{H}^{-1}(\cdot)\Phi_t(\cdot)\hat{T}_t^{-1}J_rz,
\]
we infer by \Cref{lem:funmat} and the definition of $\hat{W}$ (also use \Cref{eq:LGH}) for all $z\in \C^r$ and such $t$
\begin{equation}\label{eq:AmaxR}
A_{\max} R_tz = \i t R_t z\text{ and }
(0,z) = J_r z =   \hat{W} (\gamma  R_tz)= (L_\gamma R_tz,G_\gamma R_tz).
\end{equation}
Finally, for such $t$, we set
\[
   \Psi_t \coloneqq \int_a^b \Phi_t(x)^*\mathcal{H}(x)^{-1}\Phi_t(x)\d x\text{ and } S_t\colon \C^r\ni  z \mapsto H_\gamma (R_t z).
\]
For the techniques provided in the next two results, we refer particularly to \cite[Corollary 4.10 and Lemma 4.11]{BMNW20}.
\begin{lemma} Assume $\i\R\subseteq \rho(A)$, $s,t\in \R$. Then
\[
   1_r - S_t^*S_s = \i (s-t) R_t^*R_s; \quad S_t^*=S_t^{-1}.
\]
Moreover,
\[
   \i  S_t^*(s\mapsto S_s)'(t) = R_t^*R_t\text{ and } \|S_t'\|=\|R_t\|^2
\]
\end{lemma}
\begin{proof}
Let $y,z \in \C^r$. Then we compute using $L_\gamma R_t z =0$ and $z=G_\gamma R_t z$ (see \Cref{eq:AmaxR})
\begin{align*}
  \i (s-t) \langle R_s y , R_t z \rangle 
&  = \langle \i s R_s y , R_t z \rangle+\langle  R_s y , \i t R_t z \rangle \\
& =\langle A_{\max}R_s y, R_t z\rangle +
\langle  R_s y , A_{\max}R_t z \rangle \\
& = \langle \Sigma \gamma R_s y, \gamma R_t z\rangle = \langle G_\gamma R_s y, G_\gamma  R_t z\rangle - \langle H_\gamma R_s y, H_\gamma  R_t z\rangle \\
& = \langle y,z\rangle-\langle S_s y,S_t z\rangle,
  \end{align*}
  which establishes the first formula. Taking $t=s$, we infer that $S_t$ is unitary.
  
  Next, note that $s\mapsto S_s$ is differentiable by differentiable dependence on parameters for solutions of ordinary differential equations. Thus, we compute
  \begin{align*} S_t^*S_t' & =
       (s\mapsto S_t^*S_s)'(t) \\ &=    -  (s\mapsto(1_r - S_t^*S_s))'(t) =   -(s\mapsto(\i (s-t) R_t^*R_s))'(t) \\ &= -\i R_t^*R_t - \i (t-t)R_t^*(s\mapsto R_s)'(t)= -\i R_t^*R_t.
  \end{align*}
  Unitarity of $S_t$ shows the remaining formula.
\end{proof}
\begin{lemma}\label{lem:GP2} Assume $\i\R\subseteq \rho(A)$, let $t\in \R$. Then
\[
    (\i t -A)^{-1} + \big((\i t -A)^{-1}\big)^*=R_tR_t^*
\]
and
\[
   \|R_t^*R_t\| =   \|R_tR_t^*\|\leq 2\|(\i t -A)^{-1} \|.
\]
\end{lemma}
\begin{proof}
  Let $u\coloneqq (\i t -A)^{-1} f$ and $z\in \C^r$. Then we apply the integration by parts formula to $w\coloneqq R_t z$ and obtain
  \begin{align*}
     \langle f, R_t z\rangle & = \langle (\i t -A) u, R_t z\rangle \\
     & =\langle \i t u, R_t z\rangle - \langle A u, R_t z\rangle \\
     & =\langle  u, -A_{\max}R_t z\rangle - \langle A_{\max} u, R_t z\rangle \\ 
     & = -\langle G_\gamma u,G_\gamma R_t z \rangle + \langle H_\gamma u, H_\gamma R_tz\rangle \\
     & = \langle H_\gamma  u, S_t z\rangle\\
     & = \langle H_\gamma  (\i t -A)^{-1} f, S_t z\rangle.
  \end{align*}
  Thus, 
  \[
     R_t^* f = S_t^* H_\gamma(\i t -A)^{-1} f
  \]and, by unitarity of $S_t$,
  \[
  S_t R_t^* f = H_\gamma(\i t -A)^{-1} f.
  \]
  Next,
  \[
    2\Re \langle f, (\i t -A)^{-1} f\rangle = -2\Re \langle A u,u\rangle = \|H_\gamma (\i t -A)^{-1} f\|^2 = \|S_tR_t^*f\|^2= \|R_t^*f\|^2.
  \]
  Hence, 
  \[
  (\i t -A)^{-1}+ \big((\i t -A)^{-1}\big)^*= R_t R_t^*.
  \]
  Finally, the last equality follows from the previous.
\end{proof}

\begin{proposition} Assume $\i\R\subseteq \rho(A)$ and let $t\in \R$.

Then
\[
    s_t \coloneqq \|(\i t -A)^{-1}\|^{-1}= \min\{s>0; \det (S_{t+s}+S_{t-s})=0\}.
\]
Moreover,
\[
    \int_{t-s_t}^{t+s_t} \|R_s^*R_s\| \d s\geq 2.
\]
\end{proposition}
\begin{proof}
By compactness of $(\i t-A)^{-1}$, we find a unit vector $w\in X$ such that $\|(\i t-A)^{-1}w\| = s_t^{-1}$. It follows that
\[
s_t^{-2}\langle w,w\rangle = s_t^{-2} = \|(\i t-A)^{-1}w\|^2 = \langle \big((\i t-A)^{-1}\big)^*(\i t-A)^{-1} w,w\rangle
\] and, hence 
\[
\big((\i t-A)^{-1}\big)^*(\i t-A)^{-1} w = s_t^{-2}w.
\] With $u\coloneqq s_t (\i t-A)^{-1} w$ we get 
\[
   (\i t-A) u = s_t w.
\]Moreover,
\[
  s_t ( (\i t-A)^{-1})^* u =  s_t ( (\i t-A)^{-1})^* s_t (\i t-A)^{-1} w= w; \text{ i.e., }(\i t-A)^{*} w = s_t u.
\]We note in passing that, in particular, $w \in \dom(A^*)$ and $u\in \dom(A)$.

The above equations read
\[
A_{\max} u = \i t u - s_t w\text{ and } A_{\max}w= \i t w +s_t u.
\]
 Next, we put $y_\pm \coloneqq u\pm \i w$. Then
\[
A_{\max} y_+= \i t u - s_t w + \i (\i t w+s_t u) = \i (t+s_t) (u+\i w) = \i (t+s_t)y_+
\]
and, similarly,
\[
A_{\max} y_-= \i t u - s_t w - \i (\i t w+s_t u) =   \i (t-s_t)y_-.
\]
Since $w,u \in \dom(A^*)\cup \dom(A)$, we infer that $L_\gamma y_\pm =0$. As $G_\gamma u=0$ and $H_\gamma w=0$ it follows that
\[
G_\gamma y_+=-G_\gamma y_-\text{ and }H_\gamma y_+=H_\gamma y_-.
\]
Next, we define $z\coloneqq G_\gamma y_+$. Then $z\neq0$. Indeed, $z=0$ implies $w\in \dom(A)$ and, thus, $y_\pm \in \dom(A)$. Next, since $w\neq 0$,  $y_\pm$ cannot both be simultaneously zero; thus one of them is an eigenvector on $\i\R$, contradicting the assumption.

Then we compute
\[
  A_{\textnormal{max}} R_{t+s_t} z = \i (t+s_t)R_{t+s_t} z\text{ and }\hat{W}\gamma R_{t+s_t} z = (0,z). 
\]
Thus, by uniqueness of the (1D) boundary value problem, we have $R_{t+s_t}z=y_+$. Similarly, $\hat{W}\gamma R_{t-s_t}z= (0,z)$ leading to $R_{t-s_t}z =-y_-$. Hence, 
\[
   \big(S_{t+s_t}+S_{t-s_t}\big)z = H_\gamma  (R_{t+s_t} z+R_{t-s_t}z) = 0.
\]
This proves that $s_t \in \{s>0; \det(S_{t+s}+S_{t-s})=0\}$.

Next, let $s>0$ and $z\neq 0$ a unit vector satisfy $(S_{t+s}+S_{t-s})z=0$. Then we let
\[
   y_+ \coloneqq R_{t+s} z,\quad y_-\coloneqq -R_{t-s}z,\quad u=\frac{1}{2}(y_++y_-),\quad w=\frac{1}{2\i }(y_+-y_-).
\]
It is not difficult to see that $G_\gamma u=0$ and $H_\gamma w=0$ leading to $u\in \dom(A)$ and $w\in \dom(A^*)$. Moreover, it is not difficult to see that
\[
   (\i t -A)u = s w\text{ and }(\i t -A)^*w = s u.
\]
Again, by construction, not both $u$ and $w$ can be zero. Hence, by the invertibility of $(\i t-A)$ (and its adjoint), neither of them is zero. It follows that $ w = \big((\i t -A)^*\big)^{-1} su =s^2  \big((\i t -A)^*\big)^{-1}(\i t -A)^{-1} w$. Hence, $s^{-2}\leq \|\big((\i t -A)^*\big)^{-1}\|^2= s_t^{-2}$ and $s_t \leq s$ proving the minimum formula.

Finally, choose a unit vector $z\in \ker((S_{t+s_t}+S_{t-s_t}))$. Using the unitarity of $S_t$, we compute
\begin{align*}
 2 & = 2\|S_{t+s_t} z\|\\
 & =\|(S_{t+s_t} - S_{t-s_t})z+(S_{t+s_t}+S_{t-s_t})z\| \\
 & =\|(S_{t+s_t} - S_{t-s_t})z\| \\
 & = \big\| \int_{t-s_t}^{t+s_t} S'_\tau z \d \tau\big \| \\
 &\leq \int_{t-s_t}^{t+s_t} \|S'_\tau \| \d \tau = \int_{t-s_t}^{t+s_t} \|R_\tau\|^2 \d \tau.\qedhere
\end{align*}
\end{proof}
\begin{remark}
The lower bound of $2$ can be improved to $\pi$ as $t\mapsto S_t$ is a continuously differentiable path of unitaries; see also \cite{ALV14}. We furthermore refer to \cite[Prop. I.2.1 and I.1.20(1)]{BH99}, together with 
\cite[Ex. I.3.6(1)]{BH99}, to infer that every curve on $\mathbb{S}^{2r-1}\subseteq \R^{2r}$ joining antipodal points has length at least $\pi$.

For convenience, we include the proof in the present situation. For a unit vector $z\in\ker(S_{t+s_t}+S_{t-s_t})$ we define $\eta\colon[t-s_t,t+s_t]\ni \tau\mapsto S_{\tau}z$. By unitarity, $\|\eta(\tau)\|^2=1$. In particular, differentiation w.r.t~$\tau$ leads to $\Re\langle\eta'(\tau),\eta(\tau)\rangle=0$. Hence, for $\eta_- \coloneqq \eta(t-s_t)$ and $f(\tau)=\Re\langle\eta(\tau),\eta_-\rangle$, we get
\[
 |f'(\tau)|
 =\big|\Re\langle\eta'(\tau),\eta_--f(\tau)\eta(\tau)\rangle\big|
 \leq \|\eta'(\tau)\|\sqrt{1-f(\tau)^2}.
\]
For $0<\varepsilon<1$, the function $\theta_\varepsilon=\arccos((1-\varepsilon)f)$ therefore satisfies
$|\theta_\varepsilon'|\leq\|\eta'\|$.
Since $f(t-s_t)=1$ and $f(t+s_t)=-1$, integration and then $\varepsilon\to0+$ yield
\[
 \pi\leq\int_{t-s_t}^{t+s_t}\|\eta'(\tau)\|\,\d\tau
 \leq\int_{t-s_t}^{t+s_t}\|S_\tau'\|\,\d\tau
 =\int_{t-s_t}^{t+s_t}\|R_\tau\|^2\,\d\tau.
\]
Consequently, the estimate \eqref{eq:GP1} below can be sharpened to
\[
 \sup_{t\in\R}\|(\i t-A)^{-1}\|
 \leq\frac{2}{\pi}\sup_{s\in\R}\|R_s\|^2.
\]
This refinement is not needed for the characterisation.
\end{remark}
We will prove \Cref{thm:mr-raw} in the following form. We will use 
\[
  \Psi_t \coloneqq \int_{a}^b \Phi_t(x)^*\mathcal{H}(x)^{-1} \Phi_t(x)\d x.
\]
\begin{theorem}[Characterisation of exponential stability without (B)]\label{thm:charwoB} Let $A$ be as above. Then the following conditions are equivalent:
\begin{enumerate}
\item[(i)] $A$ generates an exponentially stable semi-group;
\item[(ii)] $T_t$ is invertible for all $t\in \R$ and 
\[
  \sup_{t\in \R} \|R_t\|<\infty.
\]
\item[(iii)]$T_t$ is invertible for all $t\in \R$ and 
\[
  \sup_{t\in \R} \|\Psi_t^{1/2} \hat{T}_t^{-1} J_r \|<\infty,
\]
\end{enumerate}
\end{theorem}
\begin{proof}
Note that invertibility of $T_t$ is necessary for all (i), (ii), (iii) as it is trivially implied by either (ii) or (iii) and using \Cref{thm:WZ}, the same property also follows from (i). Moreover, by the same argument, we infer that $\i\R\subseteq \rho(A)$. 

Using the Gearhart--Pr\"uss theorem, we are left with showing that the conditions (ii) and (iii) are equivalent to uniform boundedness of the resolvent. For this let $t\in \R$ and we see that
\[
  2\leq \int_{t-s_t}^{t+s_t} \|R_s^*R_s\|\d s \leq 2s_t\sup_{s\in \R}  \|R_s\|^2.
\]
Thus, as $s_t = \|(\i t-A)^{-1}\|^{-1}$
\[
\| (\i t- A)^{-1}\|\leq \sup_{s\in \R}\|R_s\|^2.
\]
Computing the supremum over $t\in \R$, we obtain
\begin{equation}\label{eq:GP1}
  \sup_{t\in \R} \| (\i t- A)^{-1}\|\leq \sup_{s\in \R}\|R_s\|^2.
\end{equation}
It follows directly from \Cref{lem:GP2} that
\begin{equation}\label{eq:GP2}
\sup_{s\in \R}\|R_s\|^2\leq 2   \sup_{t\in \R} \| (\i t- A)^{-1}\|.
\end{equation}
Hence, (i)$\iff$(ii) has been established by \Cref{eq:GP1} and \Cref{eq:GP2}.

For (ii)$\iff$(iii), we reformulate the respective expressions. Let $t\in \R$. Then for $z\in \C^r$ we estimate
\begin{align*}
  \|R_t z\|_{X}^2 & = \int_{a}^b \langle R_t(x)z, \mathcal{H}(x)R_t(x) z \rangle_{\C^d}\d x \\
  & = \int_{a}^b \langle \mathcal{H}^{-1}(x)\Phi_t(x)\hat{T}_t^{-1}J_rz, \mathcal{H}(x)\mathcal{H}^{-1}(x)\Phi_t(x)\hat{T}_t^{-1}J_rz\rangle \d x \\
    & = \int_{a}^b \langle \Phi_t(x)^*\mathcal{H}^{-1}(x)\Phi_t(x)\hat{T}_t^{-1}J_rz, \hat{T}_t^{-1}J_rz\rangle \d x \\
    & =  \langle \Big(\int_{a}^b\Phi_t(x)^*\mathcal{H}^{-1}(x)\Phi_t(x)\d x\Big)\hat{T}_t^{-1}J_rz, \hat{T}_t^{-1}J_rz\rangle \\
    & = \| \Psi_t^{1/2} \hat{T}_t^{-1} J_r z\|_{\C^d}^2.
\end{align*}
Hence,
\[
   \|R_t \| = \| \Psi_t^{1/2} \hat{T}_t^{-1} J_r\|_{\C^d}.\qedhere
\]
\end{proof}
\begin{corollary}[Characterisation of exponential stability without (B)]\label{cor:charB} Let $A$ be as above and assume $\Re \langle Au,u\rangle<0$ for all $u\in \dom(A)$ with $\gamma u\neq 0$. Then the following conditions are equivalent:
\begin{enumerate}
\item[(i)] $A$ generates an exponentially stable semi-group;
\item[(ii)] $T_t$ is invertible for all $t\in \R$ and 
\[
  \sup_{t\in \R} \|R_t\|<\infty.
\]
\item[(iii)]$T_t$ is invertible for all $t\in \R$ and 
\[
  \sup_{t\in \R} \|\Psi_t^{1/2} {T}_t^{-1}\|<\infty.
\]
\end{enumerate}
\end{corollary}
\begin{proof}
The assumption on $A$ implies $r=d$ and, hence, $J_r=I_d$. Thus, 
\[
 \| \Psi_t^{1/2} \hat{T}_t^{-1} J_r \| = \| \Psi_t^{1/2} \hat{T}_t^{-1}  \|=\| \Psi_t^{1/2} {T}_t^{-1} Y^{-1} \| \sim \| \Psi_t^{1/2} {T}_t^{-1}\|,
 \]
 with constants independent of $t$ as $Y$ is bijective.
\end{proof}

Finally, we may go back to \Cref{thm:charwoB} under the additional requirement of 
\begin{equation}\label{eq:Bprime}\tag{B'}
\sup_{t\in \R} \|\Psi_t\|<\infty
\end{equation} instead. We will demonstrate in \Cref{sec:ex} that this condition is strictly weaker than \eqref{eq:(B)}.

\begin{theorem}[A characterisation with \eqref{eq:Bprime}]\label{thm:charBpr} Let $A$ be as above and assume \eqref{eq:Bprime}. Then the following conditions are equivalent.
\begin{enumerate}
\item[(i)] $A$ generates an exponentially stable semigroup;
\item[(ii)] $T_t$ is invertible for all $t\in \R$ and 
\[
\sup_{t\in \R}\| \hat{T}_t^{-1}J_r\|<\infty.
\]
\end{enumerate}
\end{theorem}
\begin{proof}
(i)$\Rightarrow$(ii): By \Cref{thm:charwoB}, $\sup_{t\in \R}\|\Psi_t^{1/2}\hat{T}_t^{-1}J_r\|<\infty$. Next, by \cite[Lemma 3.2]{TW26}, we have $\Phi_t(x)^{-1}=P_1^{-1}\Phi_t(x)^*P_1$ for almost all $x\in (a,b)$ and, thus,
\begin{align*}
   \big(\Phi_t^*(x)\mathcal{H}(x)^{-1}\Phi_t(x)\big)^{-1} & = \Phi_t(x)^{-1}\mathcal{H}(x)(\Phi_t(x)^*)^{-1} \\
& = P_1^{-1}\Phi_t(x)^*P_1\mathcal{H}(x)P_1\Phi_t(x)P_1^{-1} \\   
& \leq  \|P_1\|^2 \|\mathcal{H}\|_{\infty}P_1^{-1}\Phi_t(x)^*\Phi_t(x)P_1^{-1} \\
& \leq  \|P_1\|^2 \|\mathcal{H}\|_{\infty}^2 P_1^{-1}\Phi_t(x)^*\mathcal{H}(x)^{-1}\Phi_t(x)P_1^{-1}.
\end{align*}
Then, by Jensen's inequality (see \Cref{lem:pdm} below; for this note that $\Psi_t$ is invertible as $\mathcal{H}$ is strictly positive definite and $\Phi_t$ is close to $1$ at $a$),
\begin{align*}
\Psi_t^{-1} & \leq \frac{1}{(b-a)^2} \int_a^b  \big(\Phi_t^*(x)\mathcal{H}(x)^{-1}\Phi_t(x)\big)^{-1} \d x \\
& \leq \frac{1}{(b-a)^2} \int_a^b  \|P_1\|^2 \|\mathcal{H}\|_{\infty}^2 P_1^{-1}\Phi_t(x)^*\mathcal{H}(x)^{-1}\Phi_t(x)P_1^{-1} \d x \\
& \leq \frac{1}{(b-a)^2}\|P_1\|^2 \|\mathcal{H}\|_{\infty}^2 P_1^{-1} \int_a^b  \Phi_t(x)^*\mathcal{H}(x)^{-1}\Phi_t(x)\d x P_1^{-1}.
\end{align*}
Hence,
\[
 \|\Psi_t^{-1}\|\leq \frac{1}{(b-a)^2}\|P_1\|^2 \|P_1^{-1}\|^2  \|\mathcal{H}\|_{\infty}^2 \|\Psi_t\|.
\]
Thus, we estimate
\begin{align*}
   \sup_{t\in \R}\| \hat{T}_t^{-1}J_r\|^2 &=  \sup_{t\in \R}\|\Psi_t^{-1/2} \Psi_t^{1/2}\hat{T}_t^{-1}J_r\|^2 \\
   & \leq \sup_{t\in \R}\|\Psi_t^{-1}\| \sup_{t\in \R}\|\Psi_t^{1/2}\hat{T}_t^{-1}J_r\|^2 \\
   & \leq \frac{1}{(b-a)^2}\|P_1\|^2 \|P_1^{-1}\|^2  \|\mathcal{H}\|_{\infty}^2\sup_{t\in \R}  \|\Psi_t\| \sup_{t\in \R}\|\Psi_t^{1/2}\hat{T}_t^{-1}J_r\|^2 <\infty.
\end{align*}
(ii)$\Rightarrow$(i) The provided estimate together with the assumption clearly implies \Cref{thm:charwoB}(iii) to eventually yield (i).
\end{proof}
We conclude this section with Jensen's inequality for matrices. For convenience of the reader, we provide the elementary proof, as the (conveniently available) results in this direction usually assume continuity instead of measurability; see \cite{HP81,HP03}.
\begin{lemma}[Jensen's inequality]\label{lem:pdm} Let $M\colon (a,b)\to \C^{d\times d}$ take values in self-adjoint, positive definite and invertible matrices and assume that $M,M^{-1}\in L_1(a,b)^{d\times d}$ and that $\int_a^bM(x)\d x$ is injective. Then
\[
   (\int_{a}^b M(x) \d x)^{-1} \leq \frac{1}{(b-a)^2} \int_a^b M(x)^{-1} \d x.
\]
In particular,
\[
  \big\|(\int_{a}^b M(x) \d x)^{-1}\big\|\leq \big\|\frac{1}{(b-a)^2} \int_a^b M(x)^{-1} \d x\big\|.
\]
\end{lemma}
\begin{proof}
Before we prove the claimed inequality, we recall that for any $C\in \C^{d\times d}$ self-adjoint, positive definite and invertible and $u\in \C^d$ we have
\[
   \langle u, C^{-1} u\rangle = \sup_{y\in \C^{d}} \big(2\Re \langle u,y\rangle - \langle y, Cy\rangle\big).
\]
Indeed, from
\[
  0\leq \langle y - C^{-1} u, C(y-C^{-1} u)\rangle = \langle y, Cy\rangle -2\Re \langle y, u\rangle + \langle u, C^{-1} u\rangle
\]
we get $\langle u, C^{-1} u\rangle\geq 2\Re \langle u,y\rangle - \langle y, Cy\rangle$ with equality at $y=C^{-1} u$.

To prove the claimed inequality, we set $B\coloneqq \frac{1}{b-a}\int_{a}^b M(x) \d x$.  Then, by the preliminary observation, for $u,y \in \C^d$, we have
\begin{align*}
2\Re \langle u,y\rangle - \langle y, B y \rangle & =\frac{1}{b-a}  \int_{a}^b \Big(2\Re \langle u,y\rangle - \langle y, M(x) y \rangle\Big) \d x \\
& \leq \frac{1}{b-a}  \int_{a}^b \sup_{y\in \C^d}\Big(2\Re \langle u,y\rangle - \langle y, M(x) y \rangle\Big) \d x \\
& =\frac{1}{b-a}  \int_{a}^b \langle u, M(x)^{-1} u \rangle \d x \\
      & = \langle u, \big( \frac{1}{b-a}  \int_{a}^b M(x)^{-1}  \d x\big)u \rangle.
\end{align*}
Thus, computing the supremum over $y\in \C^d$ on the left, we get
\[
   \langle u, B^{-1}u \rangle \leq \langle u, \big( \frac{1}{b-a}  \int_{a}^b M(x)^{-1}  \d x\big)u \rangle.
\]
Since $B^{-1}=(b-a)(\int_a^b M(x)\,\d x)^{-1}$, the claim follows.
\end{proof}

\Cref{thm:charBpr} enables us to weaken a sufficient criterion for exponential stability from \cite{TW26} with -- given the above results -- an obvious proof.

\begin{corollary}\label{cor:suffco} Let $A$ be as above and assume that \eqref{eq:Bprime} holds. If $T_t$ is invertible for all $t\in \R$ and $\sup_{t\in \R}\|T_t^{-1}\|<\infty$, then $A$ generates an exponentially stable semigroup.
\end{corollary}
\section{Examples}\label{sec:ex}

In the following, we treat the same model problem with different Hamiltonians $\mathcal{H}$ to be determined later. Let $(a,b)=(0,1)$, $d=2$,
\[
   P_1 = \begin{pmatrix} 0 & 1 \\ 1 & 0 \end{pmatrix},\quad P_0 =0.
\]
Moreover, 
\[
   W = \begin{pmatrix} 1 & 1 & 0 & 0 \\ 0 & 0 & 1 & -1 \end{pmatrix}.
\]
Then $W$ has full rank. Moreover, for $u\in \dom(A)$ and $v=\mathcal{H}u$ we obtain
\begin{align*}
   2\Re\langle Au,u\rangle_X & = \langle P_1\partial v,v\rangle + \langle v,P_1\partial v\rangle \\
   & = \langle P_1 v(1),v(1)\rangle - \langle P_1 v(0),v(0)\rangle\\
   &  =  \langle  \begin{pmatrix} 0 & 1 \\ 1 & 0 \end{pmatrix} \begin{pmatrix}v_1(1)\\ v_2(1)\end{pmatrix},\begin{pmatrix}v_1(1)\\ v_2(1)\end{pmatrix}\rangle - \langle \begin{pmatrix} 0 & 1 \\ 1 & 0 \end{pmatrix} \begin{pmatrix}v_1(0)\\ v_2(0)\end{pmatrix},\begin{pmatrix}v_1(0)\\ v_2(0)\end{pmatrix}\rangle \\
&  =  \langle  \begin{pmatrix} 0 & 1 \\ 1 & 0 \end{pmatrix} \begin{pmatrix}v_1(1)\\ - v_1(1)\end{pmatrix},\begin{pmatrix}v_1(1)\\ -v_1(1)\end{pmatrix}\rangle - \langle \begin{pmatrix} 0 & 1 \\ 1 & 0 \end{pmatrix} \begin{pmatrix}v_1(0)\\ v_1(0)\end{pmatrix},\begin{pmatrix}v_1(0)\\ v_1(0)\end{pmatrix}\rangle \\
&  = - 2|v_1(1)|^2 - 2|v_1(0)|^2 = -|v(1)|^2-|v(0)|^2\leq0
\end{align*}
with strict inequality whenever $\gamma u\neq0$. Hence $d=r=2$ and $A$ is $m$-dissipative; thus generating a strongly continuous semi-group of contractions.

Moreover, $\Phi_t$ is the fundamental solution of
\[
     v' = \i t  \begin{pmatrix} 0 & 1 \\ 1 & 0 \end{pmatrix} \mathcal{H}^{-1}v\text{ subject to }\Phi_t (0)=1_2;
\]and
\[
   T_t = W\begin{pmatrix} \Phi_t(1) \\ 1_2\end{pmatrix}.
\]
\begin{proposition}\label{prop:semi-group}
If $\mathcal{H}\in L_\infty(a,b)^{d\times d}$ satisfying $\mathcal{H}(x)=\mathcal{H}(x)^*\geq c$ for a.e.~$x\in (0,1)$ and some $c>0$, then for all $t\in \R$, 
\[
    T_t^*T_t = 1_2 +\Phi_t(1)^*\Phi_t(1), \quad \|T_t^{-1}\|\leq 1.
\]
In particular,  $\i\R\subseteq \rho(A_{\mathcal{H}})$ and $A_{\mathcal{H}}$ generates a strongly stable $C_0$-semigroup.
\end{proposition}
\begin{proof}Since $\Phi_t(x)^*P_1\Phi_t(x)=P_1$ by \cite[Lemma 3.2]{TW26}, for $t\in \R$ we compute
\begin{align*}
  T_t^*T_t & = \big(W\begin{pmatrix} \Phi_t(1) \\ 1_2\end{pmatrix}\big)^*W\begin{pmatrix} \Phi_t(1) \\ 1_2\end{pmatrix} \\
  & = \begin{pmatrix} \Phi_t(1)^* & 1_2\end{pmatrix}W^*W\begin{pmatrix} \Phi_t(1) \\ 1_2\end{pmatrix}\\
  & = \begin{pmatrix} \Phi_t(1)^* & 1_2\end{pmatrix}\begin{pmatrix} 1 & 0 \\ 1&  0 \\  0 &1  \\ 0 & -1  \end{pmatrix} \begin{pmatrix} 1 & 1 & 0 & 0 \\ 0 & 0 & 1 & -1 \end{pmatrix}\begin{pmatrix} \Phi_t(1) \\ 1_2\end{pmatrix} \\
  & = \begin{pmatrix} \Phi_t(1)^* & 1_2\end{pmatrix}
  \begin{pmatrix} 1 & 1 & 0 & 0 \\ 1 & 1 & 0 & 0 \\0 & 0 & 1 & -1 \\0 & 0 & -1 & 1\end{pmatrix}
  \begin{pmatrix} \Phi_t(1) \\ 1_2\end{pmatrix} \\
  & = \Phi_t(1)^*(1_2+P_1)\Phi_t(1) + (1_2 - P_1) \\
  & = \Phi_t(1)^*\Phi_t(1) +1_2 \geq 1_2,
\end{align*}
which proves invertibility of $T_t$ as well as $\|T_t^{-1}\|\leq 1$. Hence, $\i\R\subseteq\rho(A_{\mathcal{H}})$ and strong stability follows from \cite{WZ24}; see also \Cref{prop:TWZ} and \Cref{thm:WZ}.
\end{proof}
For the next result, we recall
\[
   \Psi_t = \int_{a}^b \Phi_t(x)^*\mathcal{H}(x)^{-1} \Phi_t(x)\d x.
\]
\begin{corollary}\label{cor:reform-Rayleigh} Let $t\in \R$. Then
\[
    \|\Psi^{1/2}_t T_t^{-1}\|^2  = \sup_{v\ne 0, v' = \i t P_1\mathcal{H}^{-1}v} \frac{\int_0^1 \langle v(x), \mathcal{H}(x)^{-1} v(x)\rangle \d x}{|v(0)|^2+|v(1)|^2}
     \] 
\end{corollary}
\begin{proof}
Since $T_t$ is invertible, we may compute with $T_t^*T_t=1_2+\Phi_t^*(1)\Phi_t(1)$
  \begin{align*}
&    \|\Psi^{1/2}_t T_t^{-1}\|^2 \\ & = \|\big(T_t^{-1}\big)^* \Psi_t T_t^{-1}\| \\
    & = \sup_{c\in \C^2, c\neq 0} \frac{1}{|c|^2} \langle \Psi_t T_t^{-1}c, T_t^{-1}c\rangle \\
    & = \sup_{c\in \C^2, c\neq 0} \frac{1}{|T_t c|^2} \langle \Psi_t c, c\rangle \\
    & = \sup_{c\in \C^2, c\neq 0} \frac{1}{\langle T_t^*T_t c,c\rangle}  \int_{a}^b \langle \Phi_t(x)^*\mathcal{H}(x)^{-1} \Phi_t(x)c,c\rangle\,\d x \\
    & = \sup_{c\in \C^2, c\neq 0} \frac{1}{\langle (1_2+\Phi_t^*(1)\Phi_t(1))c,c\rangle}  \int_{a}^b \langle \mathcal{H}(x)^{-1} \Phi_t(x) c, \Phi_t(x) c\rangle  \d x \\
        & = \sup_{c\in \C^2, c\neq 0} \frac{1}{|c|^2+|\Phi_t(1)c|^2} \int_{a}^b \langle \mathcal{H}(x)^{-1} \Phi_t(x) c, \Phi_t(x) c\rangle  \d x \\
     & = \sup_{v\ne 0, v' = \i t P_1\mathcal{H}^{-1}v} \frac{\int_0^1 \langle v(x), \mathcal{H}(x)^{-1} v(x)\rangle \d x}{|v(0)|^2+|v(1)|^2}.\qedhere
     \end{align*}
\end{proof}

\subsection{Hamiltonian densities failing (B)}\label{subsec:failure-of-(B)}

Let $R =\diag(\alpha,\alpha^{-1})$ for some $\alpha>1$, and $\mathcal{H}$ be piecewise constant on right-open, half-closed intervals, equal to $R^{-1}$ and $R$ (starting with $R^{-1}$ at $x=0$), respectively. 

We have 
\[
P_1 R P_1 R =\begin{pmatrix} 0 & \alpha^{-1} \\ \alpha & 0 \end{pmatrix} \begin{pmatrix} 0 & \alpha^{-1} \\ \alpha & 0 \end{pmatrix} = 1_2
\]
Hence,
\[
     \e^{\i t P_1 R} = \cos(t)1_2+ \i \sin(t)P_1 R.
\]
Let $(x_n)_n$ in $[0,1]$ be the sequence of jumps from $R^{-1}$ to $R$ (and vice versa). 

For illustration purposes, assume there is a single jump at $x^*\in(0,1)$. Then
\[
   \Phi_t (x) = \1_{[0,x^*]} \e^{\i t P_1 R x} + \1_{(x^*,1]}  \e^{\i t P_1 R^{-1} (x-x^*)} \e^{\i t P_1 R x^{*}}
\]
Similarly, assume that $\{x_n;n\in \N\}\subseteq [0,1)$ is strictly monotone and has infinitely many elements only accumulating at $1$
\footnote{It follows that $\mathcal{H}$ is not of bounded variation.}, with $x_0=0$ we have, for $0\leq x<1$,
\[
   \Phi_t (x) = \sum_{k=0}^\infty \1_{[x_k,x_{k+1})}  \e^{\i t P_1 R^{(-1)^k} (x-x_k)}F_{k-1}(t) \cdots F_0(t),
   \]
   where, for $k\in \Nnull$,
   \[
F_k(t) \coloneqq \e^{\i t P_1 R^{(-1)^k } d_{k}} \text{ and }d_k\coloneqq x_{k+1}-x_k>0.
\]
The empty product is $1_2$. Then, for all $k\in \Nnull$,
\[
   \Phi_t (x_{k+1})= F_k(t)F_{k-1}(t) \cdots F_0(t)\text{ and }
    \Phi_t (1)= \lim_{k\to\infty}  F_k(t)F_{k-1}(t) \cdots F_0(t).
\]With these formulas, we are equipped to establish a violation of \eqref{eq:(B)} in the present setting.
\begin{proposition}\label{prop:phit1} Assume there exists an odd integer $p\geq3$ such that $d_k/d_{k+1}=p$ for all $k\in\Nnull$. Then there exists a sequence $(t_k)_k$ in $(0,\infty)$ strictly monotone increasing such that for  $k\in \N$
\[
 \|\Phi_{t_k} (x_{k+1})\| \to \infty\quad(k\to\infty) \text{ and }\sup_{t\in \R} \|\Phi_t\|_{\infty,[0,x_{k+1}]} =\|\Phi_{t_k} (x_{k+1})\|.
\]
\end{proposition}
\begin{proof}
For $k\in \N$ we define 
\[
t_k \coloneqq \frac{\pi}{2}\frac{p^{k+1}}{d_0}.
\]
Then for $0\leq \ell \leq k$ we infer using a telescoping product argument
\[
   t_k d_{\ell} =\frac{\pi}{2}\frac{p^{k+1}}{d_0} d_{\ell} =\frac{\pi}{2}p^{k+1} \prod_{m=0}^{\ell-1}\frac{ d_{m+1}}{d_m} =\frac{\pi}{2}p^{k+1} p^{-\ell} \in \frac{\pi}{2}(2\N+1).\]
   Thus,
\begin{align*}
  \| \Phi_{t_k}(x_{k+1})\|& = \|F_k(t_k)F_{k-1}(t_k) \cdots F_0(t_k)\| \\ 
  & = \| \e^{\i \frac{\pi}{2} P_1 R^{(-1)^k }}\cdots \e^{\i \frac{\pi}{2} P_1 R^{(-1)^0 }}\| =\|  P_1 R^{(-1)^k }\cdots  P_1 R^{(-1)^0 } \|
  \end{align*}
  Since
  \begin{align*}
   P_1 R^{-1} P_1 R & =\begin{pmatrix} 0 & 1\\ 1 & 0\end{pmatrix}\begin{pmatrix} \alpha^{-1} & 0\\ 0 & \alpha\end{pmatrix}\begin{pmatrix} 0 & 1\\ 1 & 0\end{pmatrix}\begin{pmatrix} \alpha & 0\\ 0 & \alpha^{-1}\end{pmatrix} =\begin{pmatrix}  0& \alpha\\ \alpha^{-1} & 0\end{pmatrix}\begin{pmatrix}  0& \alpha^{-1}\\ \alpha& 0\end{pmatrix} \\
   & = \begin{pmatrix} \alpha^{2} & 0\\ 0 & \alpha^{-2}\end{pmatrix},
  \end{align*}
  we obtain 
  \[
    \| \Phi_{t_k}(x_{k+1})\|=\begin{cases} \|\begin{pmatrix} \alpha^{2} & 0\\ 0 & \alpha^{-2}\end{pmatrix}^{\frac{k+1}{2}}\|, & k \text{ odd},\\
    \|P_1R\begin{pmatrix} \alpha^{2} & 0\\ 0 & \alpha^{-2}\end{pmatrix}^{\frac{k}{2}}\|, & k \text{ even.}\end{cases}=
\alpha^{k+1}.  \]
  This establishes the first claim. For the second one, as a preliminary step, we observe
  \[
     P_1 R^{1/2} = R^{-1/2}P_1.
  \]
  Hence, 
  \[
    R^{1/2}(P_1R)R^{-1/2} =     R^{1/2}P_1R^{1/2} = P_1.
  \]
  Thus,
  \[
      \e^{\i t P_1R} = R^{-1/2} \e^{\i t P_1}R^{1/2}.
  \]
  By self-adjointness of $P_1$, $\e^{\i t P_1}$ is unitary, and thus,
  \[
    \|\e^{\i t P_1R}\| \leq \|R^{-1/2}\|\|R^{1/2}\| \leq \alpha.
  \]
  Thus, using the same argument for $R^{-1}$ replacing $R$, we deduce
  \[
    \sup_{t\in \R} \|\Phi_t\|_{\infty,[0,x_{k+1}]}\leq \alpha^{k+1},  \]
    which is attained by the argument above.
\end{proof}

\begin{remark}\label{rem:formforPhi} (a) A closer look at the proof shows that $t_k = \frac{\pi}{2}\frac{p^{2k+2}}{d_0}$ can be chosen such that, for $k\geq 1$, we have
\[
     \Phi_{t_{k}}(x_{2k}) = \varepsilon_k \begin{pmatrix} \alpha^{2k}& 0 \\ 0 &\alpha^{-2k} \end{pmatrix},
\]for some $\varepsilon_k \in \{-1,1\}$.

(b) Note that the statement in \Cref{prop:phit1} does not invoke the boundary condition matrix $W$ but rather is a statement of the mutual relationship of $P_1$ and $\mathcal{H}$.
\end{remark}
\subsection{The matrix $\Psi_t$ and an exponentially stable example}

We recall $\alpha>1$ and
\begin{align*}
\Phi_t(x) & = \sum_{k=0}^\infty \1_{[x_k,x_{k+1})}  \e^{\i t P_1 R^{(-1)^k} (x-x_k)}F_{k-1}(t)\cdots F_0(t), \\
\mathcal{H}(x)^{-1}& = \sum_{k=0}^\infty \1_{[x_k,x_{k+1})}(x) R^{(-1)^k}\text{ and }\Psi_t  =\int_0^1 \Phi_t(x)^*\mathcal{H}(x)^{-1}\Phi_t(x)\d x.
\end{align*}

\begin{proposition}\label{prop:piecewise} Let $t\in \R$, $k\in \Nnull$. Then for all $x_k\leq x<x_{k+1}$ the following holds:
\begin{enumerate}
\item[(a)]
\[
  E_k(t) \coloneqq \Phi_t(x_k)^*R^{(-1)^k}\Phi_t(x_k) = \Phi_t(x)^*\mathcal{H}(x)^{-1}\Phi_t(x) 
\]
\item[(b)]
\[
   E_0(t) = R \leq \alpha 1_{2},\quad E_{k+1}(t)\leq \alpha^2 E_k(t).
\]
%\item[(c)]
%\[
% \|\Phi_t(x_{k+1})^*R^{(-1)^{k+1}}\Phi_t(x_{k+1})\|\leq \|R\|^2\|\Phi_t(x_{k})^*R^{(-1)^{k}}\Phi_t(x_{k})\|.
%\]
%and
%\[
% \|\Phi_t(x_{k})^*R^{(-1)^{k}}\Phi_t(x_{k})\|\leq \|R\|^4  \|\Phi_t(x_{k+1})^*R^{(-1)^{k+1}}\Phi_t(x_{k+1})\|.
%\]
\end{enumerate}
\end{proposition}
\begin{proof}
(a)  It is not difficult to see that for any $t\in \R$,
\[
  \big( \e^{\i t P_1 R^{\pm 1}}\big)^* =    \e^{\big(\i t P_1 R^{\pm 1}\big)^*}=
   \e^{-\i t R^{\pm 1}P_1 }=   \e^{-\i t P_1R^{\mp 1} }.
\]
Next, recalling $P_1R^{\pm1} =R^{\mp 1/2}P_1R^{\pm 1/2}$, we obtain
\[
      \e^{\i t P_1R^{\pm1}} = R^{\mp 1/2} \e^{\i t P_1}R^{\pm 1/2}.
\]
Hence, 
\begin{equation}\label{eq:e0}
 \e^{-\i t P_1R^{\mp1} }R^{\pm1}\e^{\i t P_1 R^{\pm1}} =   R^{\pm 1/2} \e^{-\i t P_1}R^{\mp 1/2}R^{\pm1} R^{\mp 1/2} \e^{\i t P_1}R^{\pm 1/2} = R^{\pm1}.
\end{equation}
As for $k\in \Nnull$ and $x_k\leq x<x_{k+1}$ we get
\[
\Phi_t(x) =  \e^{\i t P_1 R^{(-1)^k} (x-x_k)}\Phi_t(x_k)\text{ and, thus, } \Phi_t(x)^* =\Phi_t(x_k)^*\e^{-\i t P_1 R^{-(-1)^k} (x-x_k)},
\]
 we have
\begin{align*}
\Phi_t(x)^* \mathcal{H}(x)^{-1}\Phi_t(x) & = \Phi_t(x_k)^*\e^{-\i t P_1 R^{-(-1)^k} (x-x_k)}R^{(-1)^k} \e^{\i t P_1 R^{(-1)^k} (x-x_k)}\Phi_t(x_k) \\ & =\Phi_t(x_k)^*R^{(-1)^k}\Phi_t(x_k).
\end{align*}
(b) The first equality follows from $\Phi_t(0)=1_2$. For the second, we see that
\[
    R^{\pm 1} \leq \alpha^2 R^{\mp 1}
\]
and obtain, with (a), in the sense of positive definiteness, using continuity of $\Phi_t(\cdot)$,
\[
   E_{k+1}(t) = \Phi_t(x_{k+1})^*R^{(-1)^{k+1}}\Phi_t(x_{k+1}) \leq \alpha^2 \Phi_t(x_{k+1})^*R^{(-1)^{k}}\Phi_t(x_{k+1}) =\alpha^2 E_{k}(t).\qedhere
\]
%(c) We compute 
%\begin{align*}
%  \|\Phi_t(x_{k+1})^*R^{(-1)^{k+1}}\Phi_t(x_{k+1})\|& =  \|R^{\frac{1}{2}(-1)^{k+1}}\Phi_t(x_{k+1})\|^2 \leq \|R\|\|\Phi_t(x_{k+1})\|^2\\
%  & =\|R\| \|\e^{\i t P_1 R^{(-1)^k} (x_{k+1}-x_k)}\Phi_t(x_k)\|^2 \\
%  & =\|R\| \|R^{- 1/2(-1)^k} \e^{\i t(x_{k+1}-x_k) P_1}R^{1/2(-1)^k}\Phi_t(x_k)\|^2 \\
%  & \leq \|R\|^2\|R^{1/2(-1)^k}\Phi_t(x_k)\|^2\\
%  & =\|R\|^2\|\Phi_t(x_{k})^*R^{(-1)^{k}}\Phi_t(x_{k})\|.
%\end{align*}
%Moreover,
%\begin{align*}
%&  \|\Phi_t(x_{k})^*R^{(-1)^{k}}\Phi_t(x_{k})\| \\ 
%& = \|R^{1/2(-1)^k}\Phi_t(x_k)\|^2 \\
% & =\| \big(R^{-1/2(-1)^k} \e^{\i t P_1 R^{(-1)^k} (x_{k+1}-x_k)}R^{-1/2(-1)^k} \big)^{-1}  R^{1/2(-1)^{k+1}} \Phi_t(x_{k+1})\|^2  \\
% & \leq \|R\| \|R\|^2 \|R\| \|R^{1/2(-1)^{k+1}} \Phi_t(x_{k+1})\|^2.\qedhere
%\end{align*}
\end{proof}

\begin{corollary}\label{cor:exp_stable} In the setting provided in \Cref{sec:ex} and \Cref{subsec:failure-of-(B)} and corresponding $\mathcal{H}$ with $d_k/d_{k+1}=p$ for some odd $p\in \N$ and all $k\in \Nnull$ with $p>\alpha^2$, then $A_{\mathcal{H}}$ generates an exponentially stable semi-group.
\end{corollary}
\begin{proof}
By \Cref{prop:semi-group}, we infer that $A_{\mathcal{H}}$ generates a strongly stable semi-group (with $T_t$ invertible) and, in fact, $\|T_t^{-1}\|\leq 1$. Thus, it suffices to see that \Cref{cor:charB}(iii) holds. And this, in turn, can be established if $\sup_{t\in \R}\|\Psi_t\|<\infty$ is satisfied. Using \Cref{prop:piecewise} and $\sum d_k=1$ so that $d_k = c p^{-k}$ for some $c>0$, we compute for $t\in \R$,
\begin{align*}
   \Psi_t  & =  \int_0^1 \Phi_t(x)^*\mathcal{H}(x)^{-1}\Phi_t(x)\d x \\
   & = (x_1-x_0)E_0(t) + \sum_{k=1}^\infty (x_{k+1}-x_k) E_k(t) \\
   & \leq d_0 E_0(t) + \sum_{k=1}^\infty d_k \alpha^{2k} E_{0}(t) \leq \alpha \sum_{k=0}^\infty c \big(\frac{\alpha^2}{p}\big)^k<\infty,
\end{align*}
which establishes the assertion as $p>\alpha^2$.
\end{proof}

\begin{remark}
The proof of \Cref{cor:exp_stable} proves that $\sup_{t\in \R}\|\Psi_t\|<\infty$ similarly to \Cref{rem:formforPhi} (b) without using the boundary condition matrix $W$. Hence, the set-up provided confirms the statement at the end of \Cref{sec:res} that, in the case for $r=d$, \Cref{thm:charBpr} is stronger than the main result in \cite{TW26}. In any case, \Cref{cor:suffco} weakens the sufficient direction of \Cref{thm:WTB} by replacing \eqref{eq:(B)} with \eqref{eq:Bprime}.
\end{remark}

\subsection{A mirrored example failing to be exponentially stable}

Let $\mathcal{H}$ be given as in \Cref{subsec:failure-of-(B)} with a sequence $(x_n)_n$ of strictly monotonically increasing jumping points only accumulating at $1$; moreover, we will assume that for $d_k=x_{k+1}-x_k$ we find $p\in \N$, $p\geq 3$ odd such that for all $k\in \Nnull$
\[
    \frac{d_k}{d_{k+1}} = p.
\] Throughout this section, we again let $\alpha>1$. Then we define 
\[
\mathcal{K} \colon x\mapsto \mathcal{H}(|2x-1|),
\]
which retains the boundedness and positive-definiteness conditions of $\mathcal{H}$. In particular, by \Cref{prop:semi-group}, $A_{\mathcal{K}}$ generates a strongly stable $C_0$-semigroup. Moreover, by the boundary conditions, it still follows that $\Re\langle A_{\mathcal{K}} u,u\rangle<0$ for all $u\in \dom(A_{\mathcal{K}})$ with $\gamma u\neq 0$, yielding applicability of \Cref{cor:charB}. In the following, we will show that (iii) fails for $A_{\mathcal{K}}$: 

\begin{proposition} The semigroup generated by $A_{\mathcal{K}}$ fails \Cref{cor:charB} (iii) and (B). In particular, the semi-group is not exponentially stable.
\end{proposition}
\begin{proof} To further distinguish the setting with $\mathcal{H}$ and $\mathcal{K}$, we shall append the corresponding index to the objects initially defined for $\mathcal{H}$, only.
By \Cref{rem:formforPhi}, we find $(t_k)_k$ in $(0,\infty)$ such that $\Phi_{\mathcal{H},t_k}(x_{2k}) =\varepsilon_k\diag(\alpha^{2k},\alpha^{-2k})$ for some $\varepsilon_k \in \{-1,1\}$. Next, define $w_k \coloneqq \Phi_{\mathcal{H},t_k}\begin{pmatrix} 0 \\ 1\end{pmatrix}$ for all $k\in \N$. For $D=\diag(-1,1)$, we define
\[
v_k\colon    [0,1]\ni x\mapsto  \begin{cases} D w_{k}(1-2x),& 0\leq x\leq \frac{1}{2},\\ w_k(2x-1),& \frac{1}{2} < x\leq 1.\end{cases}
\]
By definition $v_k$ is $H^1$ on $(0,\frac{1}{2})$ and $(\frac{1}{2},1)$ as well as continuous at $1/2$, as 
\[
w_k(0)=\Phi_{\mathcal{H},t_k}(0)\begin{pmatrix} 0 \\ 1\end{pmatrix}=\begin{pmatrix} 0 \\ 1\end{pmatrix}=D\begin{pmatrix} 0 \\ 1\end{pmatrix}
\] implying that $v_k\in H^1(0,1)$ and $v_k(1/2)=(0,1)$.

Moreover, 
\begin{align*}
   v'_k(x) & =\begin{cases}  -D 2 w_{k}'(1-2x), & 0\leq x\leq \frac{1}{2}, \\
      2w_k'(2x-1), & \frac{1}{2}<x\leq 1. \end{cases} \\
       & = 2 \i t_k
       \begin{cases}  -D  \begin{pmatrix} 0 & 1 \\ 1 & 0 \end{pmatrix}  \big(\mathcal{H}(1-2x)\big)^{-1}w_k(1-2x), & 0\leq x\leq \frac{1}{2}, \\
    \begin{pmatrix} 0 & 1 \\ 1 & 0 \end{pmatrix}    \big(\mathcal{H}(2x-1)\big)^{-1}w_k(2x-1), & \frac{1}{2}<x\leq 1. \end{cases} \\
        & = 2 \i t_k \begin{pmatrix} 0 & 1 \\ 1 & 0 \end{pmatrix}
       \begin{cases}  D   \mathcal{K}(x)^{-1}w_k(1-2x), & 0\leq x\leq \frac{1}{2}, \\
\mathcal{K}(x)^{-1}       w_k(2x-1), & \frac{1}{2}<x\leq 1. \end{cases} \\
 & = 2 \i t_k \begin{pmatrix} 0 & 1 \\ 1 & 0 \end{pmatrix}\mathcal{K}(x)^{-1}  
       v_k(x)
\end{align*}
Moreover, $v_k(0)=Dw_k(1)$ and $v_k(1)=w_k(1)$. With this, we use \Cref{cor:reform-Rayleigh} and write
\begin{align*}
    \|\Psi_{\mathcal{K},2t_k}^{1/2} T_{\mathcal{K},2t_k}^{-1}\|^2  & = \sup_{v\ne 0, v' = \i 2t_k P_1\mathcal{K}^{-1}v} \frac{\int_0^1 \langle v(x), \mathcal{K}(x)^{-1} v(x)\rangle \d x}{|v(0)|^2+|v(1)|^2} \\
    & \geq \frac{\int_0^1 \langle v_k(x), \mathcal{K}(x)^{-1} v_k(x)\rangle \d x}{|v_k(0)|^2+|v_k(1)|^2} \\
    & = \frac{\int_0^1 \langle w_k(x),\mathcal{H}(x)^{-1}w_k(x)\rangle \d x}{|Dw_k(1)|^2+|w_k(1)|^2} \\
    & \geq \frac{d_0}{2\alpha |w_k(1)|^2}. 
     \end{align*}
     Thus, for the failure of \Cref{cor:charB} (iii), it remains to show that $|w_k(1)|\to 0$. For this, by \Cref{rem:formforPhi}, we infer
     \[
         w_k(x_{2k}) = \Phi_{\mathcal{H},t_k}(x_{2k})\begin{pmatrix} 0 \\ 1 \end{pmatrix} = \varepsilon_k \alpha^{-2k}\begin{pmatrix} 0 \\ 1 \end{pmatrix}.
     \]
From $d_k/d_{k+1}=p$ it follows that $d_k = (1-p^{-1}) p^{-k} $ as $\sum_{k\in \Nnull} d_k = 1$. Thus,
    \[
    1-x_{2k} = 1 - \sum_{\ell=1}^{2k} (x_{\ell}-x_{\ell-1}) =1-\sum_{\ell=1}^{2k} d_{\ell-1} = 1- \sum_{\ell=0}^{2k-1}  (1-p^{-1}) p^{-\ell} =1 + p^{-{2k}}-1 =p^{-{2k}}.
    \]
Since $w_k' = \i t_k P_1 \mathcal{H}^{-1} w_k$, we infer from $t_k p^{-2k} = \frac{\pi}{2}\frac{p^{2}}{d_0} \eqqcolon \tau$,
    \[
       |w_k(1)| \leq \e^{t_k \|P_1\mathcal{H}^{-1}\|_{\infty}p^{-2k}}|w_{k}(x_{2k})| \leq \e^{\tau \alpha}\alpha^{-2k} \to 0\quad(k\to\infty).
    \]
    Finally, as $\Phi_{\mathcal{K},2t_k}(1/2)v_k(0) = v_k(1/2)=(0,1)$,
    \[
       \|\Phi_{\mathcal{K},2t_k}(1/2)\|\geq \frac{|v_k(1/2)|}{|v_k(0)|} = \frac{1}{|w_k(1)|} \geq \e^{-\tau \alpha}\alpha^{2k} \to \infty\quad(k\to\infty),
    \] proving the failure of (B).
\end{proof}

\section*{Data Availability}

Data sharing does not apply to this article as no datasets were generated or analysed in this study.

\section*{Conflicts of Interest}

The authors declare no conflicts of interest.

\section*{Declaration on the usage of AI tools}

The new quantity in the characterisation of exponential stability and the key estimate concerning the bound of the resolvent have been developed after communication with ChatGPT Pro 5.6. The principle of the (counter-)examples has been suggested by ChatGPT Pro 5.6.  The proofs and examples in their final form were developed in full, verified, and written by the authors, who take full responsibility for all claims and arguments in the paper.

\enlargethispage{\baselineskip}

\end{document}